\documentclass[12pt]{amsart}
\usepackage{amscd,amsthm,amsfonts,amssymb,amsmath,euscript}
\usepackage[dvips]{graphicx}
\usepackage[dvips]{graphics}
\usepackage[T2A]{fontenc}
\usepackage[cp866]{inputenc}

\newtheorem{theorem}{Theorem}[section]
\newtheorem{lemma}{Lemma}[section]

\newtheorem{note}{Remark}[section]

\newcommand{\Aut}{\mathop{\sf Aut}\nolimits}

\newcommand{\tr}{\mathop{\sf tr}\nolimits}
\newcommand{\End}{\mathop{\sf End}\nolimits}

\newcommand{\Hom}{\mathop{\sf Hom}\nolimits}

\begin{document}

\centerline {\bf Cyclic Foam Topological Field Theories}

{\ }

\centerline {\bf Sergey M. Natanzon{\footnote{Supported by grants
RFBR-07-01-00593, NWO 047.011.2004.026 (05-02-89000-NWO-a),
NSh-4719.2006.1, INTAS 05-7805}}}

{\ }

\centerline{A.N.Belozersky Institute, Moscow State University}

\centerline{Independent University of Moscow}

\centerline{Institute Theoretical and Experimental Physics}

\centerline{natanzon@mccme.ru}

{\ }

\centerline {\bf Abstract}

This paper proposes an axiomatic for Cyclic Foam Topological Field
Theories. That is Topological Field Theories, corresponding to
String Theories, where particles are arbitrary graphs. World
surfaces in this case are two-manifolds with one-dimensional
singularities. I prove that Cyclic Foam Topological Field Theories
one-to-one correspond to graph-Cardy-Frobenius algebras, that are
families $(A,B_\star,\phi)$, where $A=\{A^s|s\in S\}$ are families
of commutative associative Frobenius algebras,
$B_\star=\bigoplus_{\sigma\in\Sigma}B_\sigma$ is an graduated by
graphes, associative algebras of Frobenius type and
$\phi=\{\phi_\sigma^s: A^s\rightarrow \End(B_\sigma)|s\in
S,\sigma\in \Sigma\}$ is a family of special representations. There
are constructed examples of Cyclic Foam Topological Field Theories
and its graph-Cardy-Frobenius algebras

{\ }

\tableofcontents

\section{Introduction}\label{r0}

Two-dimensional Topological Field Theories were introduced by Segal
\cite{Se}, Atiyah \cite{At} and Witten \cite{W}. An example of them
is a topological approach of the String Theory. The String Theory is
a modern variant of uniform Field Theory. It treads particles as
one-dimensional objects. A path of a particle is represented by a
world surface, that is, a two-dimensional space. The Topological
String Theory assumes that a probability of a world surface depends
only on the state of the particle at the moments of birth/death and
on the topological type of the world surface. The standard
properties of the measure on world lines extend to properties of
Topological String Theory \cite{Du}.

The Topological Field Theory is a direct axiomatization of the
Topological String Theory. The Topological Field Theory is a
function on the set of two-dimensional spaces (world surface in
String Theory) endowed with marked points (points of birth/death in
String Theory) and also with vectors in the marked points (vectors
of states in String Theory). The function depends from the vectors
linearly. The properties of Topological String Theory can be
reformulated as properties of the function by a surgery of surfaces.

The simplest model treads particles as closed contours. Thus its
world surface is a closed surface, that is, two-dimensional
topological manifold without boundary. The corresponding Topological
Field Theory was constructed in \cite{At}, \cite{D2} for orientable
surfaces and in \cite{AN} for arbitrary (orientable and
non-orientable) surfaces.

In this case, the values of a Topological Field Theory on spheres
with one, two and three marked points determine the values of the
Topological Field Theory on all oriented surfaces. As well, the
values of a Topological Field Theory on spheres with one, two and
three marked points are structure constants for some associative,
commutative Frobenius algebra $A$ with a unit. Moreover, this
construction gives a one-to-one correspondence between Topological
Field Theories on closed orientable surfaces and associative,
commutative Frobenius algebras with unit \cite{D2}. To extend a
Topological Field Theory to non-orientable surfaces one has to add
new structures to $A$, namely, an involution of $A$ and an element
$U\in A$ which defines the value of the Topological Field Theory on
a projective sphere with a marked point \cite{AN}.

A String Theory where particles are closed contours and segments is
called an Open-Closed Theory. In this case world surfaces are
two-dimensional topological manifolds without boundary or with a
boundary consisting of closed contours. For orientable surfaces of
this type Topological Field Theory was constructed in \cite{Laz},
\cite{Moore}, where it is determined a pair of associative Frobenius
algebras with unit, connected by a special homomorphism $\phi$. The
first algebra is $A$, corresponding to closed surfaces. The second
algebra $B$, which in general is non-commutative, corresponds to
disks with marked points at the boundary. The one-to-one
correspondent between Open-Closed Topological Field Theories and the
families $(A,B,\phi)$ was proved in \cite{AN} and later,
independently, in \cite{Lauda}, \cite{MS}.

For orientable and non-orientable surfaces with boundary a
Topological Field Theory was constructed in \cite{AN}. We call it
the Klein Topological Field Theory. Klein Topological Field Theories
are in one-to-one correspondence with the Cardy-Frobenius algebras
which are the tuples $(A,B,\phi)$ with equipments \cite{AN}.

In the present paper I construct Cyclic Foam  Topological Field
Theories that correspond to String Theories where particles are
arbitrary graphs (for physical motivation see \cite{B}, and also
\cite{O} and references there in). In this case the world surface is
a CW-complex glued from finitely many surfaces ("patches") with
boundaries by gluing some segments of the boundaries. The boundaries
of surfaces form the singular part of the complex is called "seamed
graph". Complexes of this type are called "world-sheet foam" or
"seamed surfaces". They appear also in A-models \cite{R},
Landau-Ginsburg models \cite{KR2} and 3-dimensional topology
\cite{K} \cite{KR1}.

In this paper I consider a special class of  seamed surfaces which I
call cyclic foam. They satisfy the following conditions: (1) glued
boundary contours  of patches have compatible orientation; (2)
different boundary contours a patch included to different connected
components of the seamed graph. I assume also that any patch has a
"color" from a set $S$ and the closures of two patches have no
intersections if they have the same color.

We start with the definition of Cyclic Foam  Topological Field
Theories on cyclic foams where all patches are disks (Section
\ref{r1}). We prove that such Topological Field Theories are in
one-to-one correspondence with graph-Frobenius algebras which we
define in Section \ref{r2}. A graph-Frobenius algebra is an
associative algebra that satisfy all the properties of a Frobenius
algebra except for the finite dimensionality property. Instead, it
is presented as a sum of finite dimensional vector spaces
$B_\star=\bigoplus_{\sigma\in\Sigma}B_\sigma$ where $\Sigma$ is the
set of oriented colored graphs. The algebra $B$ from Cardy-Frobenius
algebras is the subalgebra of $B_\star$ that corresponds to the
segment.

In Section \ref{r3} we define Topological Field Theories for
arbitrary cyclic foams. Later (Section \ref{r4}) we prove that
Cyclic Foam  Topological Field Theories are in one-to-one
correspondence with families $(A,B_\star,\phi)$, where $A=\{A^s|s\in
S\}$ is a family of commutative associative Frobenius algebras with
units and $\phi=\{\phi_\sigma^s: A^s\rightarrow \End(B_\sigma)|s\in
S,\sigma\in \Sigma\}$ is a special family of representations in
$B_\sigma$.

Some classical topological objects satisfy the axioms of Topological
Field Theory. The algebraic description of Topological Field
Theories makes it possible to reduce some topological problems to
algebraic one. Such applications of Topological Field Theory appear
in the Theory of Links \cite{TT}, \cite{PT} and in the Theory of
Hurwitz Numbers.

The classical Hurwitz numbers are weighted  numbers of meromorphic
functions with prescribed topological types of critical values
\cite{H}. These numbers depend on topological types of surfaces and
critical values only. The classical Hurwitz numbers  and Hurwitz
numbers of regular coverings generate Topological Field Theories on
closed surfaces \cite{D1}, \cite{CNP}. Hurwitz numbers for surfaces
with boundary are defined in \cite{AN} for coverings by surfaces
with boundary, in \cite{AN1}, \cite{AN2} for coverings by seamed
surfaces, and in \cite{AN3} for regular coverings by seamed
surfaces. In these papers we proved that each type of these Hurwitz
numbers forms a Klein Topological Field Theory and we described
their Cardy-Frobenius algebras.

In Section \ref{r5} of preset paper it is constructed  examples of
Cyclic Foam  Topological Field Theories and corresponding
graph-Cardy-Frobenius algebras. These examples extend to cyclic
foams the Klein Topological Field Theories of regular Hurwitz
numbers from \cite{AN3}.

I thank A.Alekseevskii, S.Lando, and L.Rozansky for useful
discussions.

\section{Film Topological Field Theories}\label{r1}

\subsection{Film surfaces} \label{r1.1}
In this paper a \textit{graph} is a compact simplicial complex that
consist of simplexes of dimension 1 (\textit{edges}) and dimension 0
(\textit{vertices}). An edge is either a segment or a loop depending
on the topological type of its closure. A graph is said to be
\textit{regular} if all its edges are segments.

A compact CW-complex that consists of oriented cells of dimension 2
(\textit{disks}), cells of dimension 1 (\textit{edges}) and cells of
dimension 0 (\textit{vertices}) is said to be \textit{regular} if
its edges form a regular graph. Thus regular CW-complex $\Omega$ is
defined by a set $(\check{\Omega},\Delta,\varphi)$, where
$\check{\Omega}=\check{\Omega}(\Omega)$ is a set of closed oriented
disks , $\Delta=\Delta(\Omega)$ is a regular graph and
$\varphi:\partial\check{\Omega}\to\Delta$ is a gluing map, that is,
a homeomorphism on any connected component of $\partial\Omega$ and
$\varphi(\partial\check{\Omega})=\Delta$.

A system of cyclic orders on vertexes of connected components of a
regular CW-complex $\Omega$ is called a \textit{cyclic order on}
$\Omega$, if the cyclic orders on vertexes are compatible with the
orientations of the disks $\check{\Omega}(\Omega)$. A regular
CW-complex with a cyclic order is called \textit{almost cyclic
complex}.

A connected regular graph $\gamma\subset\Omega$ on a connected
almost cyclic complex $\Omega$ is called a \textit{graph-cut} if:
\begin{itemize}
\item the restriction of $\gamma$ to any disk $\omega\in\check{\Omega}$
either is empty, or forms one of edge of $\gamma$;
\item $\gamma$ divides $\Omega$ into two connected components that
splits vertices of  $\Omega$ into two nonempty groups,
compatible with the cyclic order on $\Omega$.
\end{itemize}

A almost cyclic complex $\Omega$ is called \textit{cyclic complex}
if for any compatible with the cyclic order division of vertexes
$\Omega$ there exists a graph-cut that realize it. A small
neighborhood of a vertex $q$ of a CW-complex is a cone over a
regular \textit{vertex graph} $\sigma_q$, with orientation of edges
generated by the orientation of the disks outside the neighborhood.
It is obviously that vertex graphs of cyclic complexes are
connected.

Fix a set $S$ of \textit{colors}. A  graph (respectively CW-complex)
is called \textit{colored} if a color $s(l)\in S$ is assigned to
each of its edges (respectively disks) $l$ and all the colors are
pairwise different for any connected component. A colored cyclic
complex is called a \textit{film surface}. Vertex the graph
$\sigma_q$ of a vertices $q$ of a film surface is a colored graph,
where colors of edges generated by colors of the disks. Denote by
$\Omega_b$ the set of vertices of a film surface $\Omega$.

\subsection{Topological Field Theory}\label{r1.2}
Below we assume that all vector spaces are defined over a field
$\mathbb{K}\supset\mathbb{Q}$. Let $\{X_m | m\in M\}$ be a finite
set of $n=|M|$ vector spaces $X_m$ over the field of complex numbers
$\mathbb{C}$. The action of the symmetric group $S_n$ on
$\{1,\dots,n\}$ induces its action on the sum of the vector spaces
$\left(\oplus_{\sigma} X_{\sigma(1)}\otimes\dots\otimes
X_{\sigma(n)}\right)$ where $\sigma$ runs over the bijections
$\{1,\dots, n\}\to M$, an element $s\in S_n$ takes
$X_{\sigma(1)}\otimes\dots\otimes X_{\sigma(n)}$ to
$X_{\sigma(s(1))}\otimes\dots\otimes X_{\sigma(s(n))}$. Denote by
$\otimes_{m\in M} X_m$ the subspace of all invariants of this
action. The vector space $\otimes_{m\in M} X_m$ is canonically
isomorphic to the tensor product of all $X_m$ in any fixed order;
the isomorphism is the projection of $\otimes_{m\in M} X_m$ to the
summand that is equal to the tensor product of $X_m$ in that order.

Two regular oriented colored graphs said to be \textit{isomorphic}
if there exists a homeomorphism that maps one to the other
preserving the colors and the orientations. Denote by
$\Sigma=\Sigma(S)$ the set of all isomorphism classes of connected
oriented colored graphs. The inversions of the orientations generate
the involution $*:\Sigma\rightarrow\Sigma$. Denote it by
$\sigma\mapsto\sigma^*$.

Consider a family of finite-dimensional vector spaces
$\{B_{\sigma}|\sigma\in\Sigma\}$ and a family of tensors
$\{K_{\sigma}^\otimes\in B_{\sigma}\otimes
B_{\sigma^*}|\sigma\in\Sigma\}$. Using these data, we define now a
functor $\mathcal{V}$ from the category of film surfaces to the
category of vector spaces. This functor assigns the vector space
$V_\Omega= (\otimes_{q\in\Omega_b}B_q)$ to any film surface
$\Omega$. Here $B_q$ is the copy of $B_{\sigma_q}$ that is a vector
space with a fixed isomorphism $B_q\rightarrow B_{\sigma_q}$.

\smallskip

We are going to describe all morphisms of the monoidal category
$\mathcal{S}$ of film surfaces and morphisms of the category of
vector spaces that correspond to it.

(1) {\it Isomorphism.} Let $\phi:\Omega\to \Omega'$ be a
homeomorphism of film surfaces, preserving the cyclic orders,
orientations of disks and its colors. Define
$\mathcal{V}(\phi)=\phi_*: V_{\Omega}\to V_{\Omega'}$ as linear
operator generated by the bijections
$\phi|_{\Omega_b}:\Omega_b\to\Omega'_b$.

(2) \textit{Cut.} Let $\Omega$ be a connected film surface and
$\gamma\subset\Omega$ be a graph-cut. The graph $\gamma$ is
represented by two graphs $\gamma_+$ and $\gamma_-$ on the closure
$\overline{\Omega\setminus\gamma}$ of $\Omega\setminus\gamma$ .
Contract these graphs to points $q_+=q_+[\gamma]$ and
$q_-=q_-[\gamma]$, respectively. The contraction produces a film
surface $\Omega'=\Omega[\gamma]$. Its vertices
$\Omega'=\Omega[\gamma]$ are the vertices of $\Omega$ and the points
$q_+$, $q_-$. The cyclic order, orientation and the coloring of
$\Omega$ induces an orientation and a coloring of $\Omega'$. Thus we
can assume that $\Omega'$ is a film surface and
$V_{\Omega'}=V_{\Omega}\otimes B_{q_+}\otimes B_{q_-}$. The functor
takes the morphism $\mathcal{V}(\eta)(x)=\eta_*(x)=x\otimes
K_{\sigma}^\otimes$, where $\sigma=\sigma_{q_+}=\sigma_{q_-}^*$, to
the morphism $\eta:\Omega\to\Omega'$.

(3) The tensor product in $\mathcal{S}$ defined by the disjoint
union of surfaces $\Omega'\otimes\Omega''\to \Omega'\coprod\Omega''$
induces the tensor product of vector spaces $\theta_*:
V_{\Omega'}\otimes V_{\Omega''}\to V_{\Omega'\sqcup\Omega''}$.

The functorial properties of $\mathcal{V}$ can be easily verified.

\smallskip

Fix a tuple of vector spaces and vectors
$\{B_{\sigma},K_{\sigma}^\otimes\in B_{\sigma}\otimes
B_{\sigma^*}|\sigma\in\Sigma\}$, defining the functor $\mathcal{V}$.
A family of linear forms $\mathcal {F}= \{\Phi_\Omega:V_{\Omega}\to
\mathbb{K}\}$ defined for all film surfaces $\Omega\in\mathcal{S}$
is called a \textit{Film Topological Field Theory} if it satisfies
the following axioms:

\vskip0.6cm $1^\circ$ {\it  Topological invariance.}
$$\Phi_{\Omega'}(\phi_*(x))=\Phi_\Omega(x)$$ for any isomorphism
$\phi:\Omega\to \Omega'$ of film surfaces.

\vskip 0.6cm $2^\circ$ {\it  Non-degeneracy.}\vskip 0.3cm Let
$\Omega$ be a film surface with only two vertices $q_1$, $q_2$. Then
$\sigma_{q_2}=\sigma^*$ if $\sigma_{q_1}=\sigma$. Denote by
$(.,.)_\sigma$ the bilinear form $(.,.)_\sigma: B_\sigma\times
B_{\sigma^*}\rightarrow\mathbb{K}$, where
$(x',x'')_\sigma=\Phi_\Omega(x'_{q_1}\otimes x''_{q_2})$. Axiom
$2^\circ$ asserts that the forms $(.,.)_\sigma$ are non-degenerated
for all $\sigma\in\Sigma$.

\vskip 0.6cm $3^\circ$ {\it Cut invariance.}
$$\Phi_{\Omega'}(\eta_*(x))=\Phi_\Omega(x)$$
for any cut morphism $\eta:\Omega\to \Omega'$ of film surfaces.
\vskip 0.6cm

$4^\circ$ {\it Multiplicativity.}

$$\Phi_{\Omega}(\theta_*(x'\otimes x'))= \Phi_{\Omega'}(x')\Phi_{\Omega''}(x'')$$

for $\Omega=\Omega'\coprod\Omega''$,  $x'\in V_{\Omega'}$, $x''\in
V_{\Omega''}$.

\smallskip

Note that a Topological Field Theory defines the tensors
$\{K_{\sigma}^\otimes\in B_{\sigma}\otimes
B_{\sigma^*}|\sigma\in\Sigma\}$, since it is not difficult to prove:
\begin{lemma} \label{l1}Let $\{\Phi_\Omega\}$ be a Film Topological Field Theory.
Then $(K_\sigma^\otimes,x_1 \otimes x_2)_\sigma= (x_1,x_2)_\sigma$,
for all $x_1\in B_\sigma$, $x_2\in B_{\sigma^*}$.
\end{lemma}

\section{Graph-Frobenius algebras}\label{r2}

\subsection{Definitions} \label{r2.1}

We say that a connected film surface $\Omega$ is a
\textit{compatible surface} for colored graphes
$\sigma_1,\sigma_2,...,\sigma_n$ if these graphes are vertex graphs
of $\Omega$ and the numeration the graphs $\sigma_i$ generates the
cyclic order of vertexes of film surface $\Omega$. Denote by
$\Omega(\sigma_1,\sigma_2,...,\sigma_n)$ the set of all isomorphism
classes of compatible surfaces for $\sigma_1,\sigma_2,...,\sigma_n$.
Then $\Omega(\sigma_1,\sigma_2,...,\sigma_n)$ is either empty or
consists of a single element.

Let $\Omega(\sigma_1,\sigma_2,\sigma_3,\sigma_4)\neq\emptyset$. Then
there exist unique classes of graph-cuts
$\sigma_{(1,2|3,4)},\sigma_{(4,1|2,3)}\in\Sigma$ such that
$\Omega(\sigma_1,\sigma_2,\sigma_{(1,2|3,4)})\neq\emptyset$,
$\Omega(\sigma_{(3,4|1,2)},\sigma_3,\sigma_4)\neq\emptyset$,
$\Omega(\sigma_4,\sigma_1,\sigma_{(4,1|2,3)})\neq\emptyset$,
$\Omega(\sigma_{(4,1|2,3)},\sigma_2,\sigma_3)\neq\emptyset$ and
$\sigma_{(3,4|1,2)}=\sigma_{(1,2|3,4)}^*$,
$\sigma_{(2,3|4,1)}=\sigma_{(4,1|2,3)}^*$.

{\ }

Consider a tuple of finite dimensional vector spaces
$\{B_{\sigma}|\sigma\in\Sigma\}$. Its direct sum $B_\star =
\bigoplus_{\sigma\in\Sigma}B_{\sigma}$ is called a \textit{colored
graph-graded vector space}.

A colored graph-graded vector space with a bilinear form
$(.,.):B_\star\times B_\star\rightarrow\mathbb{K}$ and a tree-linear
form $(.,.,.):B_\star\times B_\star\times
B_\star\rightarrow\mathbb{K}$ is called a \textit{ graph-Frobenius
algebra} if
\begin{itemize}
\item $(B_{\sigma_1},B_{\sigma_2}) = 0$ for
$\sigma_1\neq\sigma_2^*$;
\item the form $(.,.)$  is not-degenerate;
\item $(B_{\sigma_1},B_{\sigma_2},B_{\sigma_3}) = 0$ for
$\Omega(\sigma_1,\sigma_2,\sigma_3)=\emptyset$
\item $\sum_{i,j}(x_1,x_2,b_i^{(1,2|3,4)})F^{ij}_{(1,2|3,4)}(b_j^{(3,4|1,2)},x_3,x_4)=
\sum_{i,j}(x_4,x_1,b_i^{(4,1|2,3)})F^{ij}_{(4,1|2,3)}(b_j^{(2,3|4,1)},x_2,x_3)$.
\end{itemize}
Here $x_k\in B_{\sigma_k}$, $\{b_i^{(s,t|k,r)}\}$ is a basis of
$B_{(s,t|k,r)}$ and $F^{ij}_{(s,t|k,r)}$ is the inverse matrix for
$(b_i^{(s,t|k,r)},b_j^{(k,r|s,t)})$.

{\ }

We will consider $B_*$ as an algebra with the multiplication
$(x_1x_2,x_3) =(x_1,x_2,x_3)$, for $x_k\in B_{\sigma_k}$. The axiom
$\sum_{i,j}(x_1,x_2,b_i^{(1,2|3,4)})F^{ij}_{(1,2|3,4)}(b_j^{(3,4|1,2)},x_3,x_4)=
\sum_{i,j}(x_4,x_1,b_i^{(4,1|2,3)})F^{ij}_{(4,1|2,3)}(b_j^{(2,3|4,1)},x_2,x_3)$
is equivalent to associativity for the algebra $B_*$. Moreover it is
a Frobenius algebra in the sense of \cite{F} if its dimension is
finite.

\subsection{One-to one correspondence}  \label{r2.2}
\begin{theorem} \label{t1} Let $\mathcal {F}=
\{\Phi_\Omega:V_{\Omega}\to \mathbb{K}\}$ be a Film Topological
Field Theory on a tuple of finite-dimensional vector spaces
$\{B_{\sigma}|\sigma\in\Sigma\}$. Then the poli-linear forms
\begin{itemize}
\item $(x',x'')=\Phi_{\Omega(\sigma_1,\sigma_2)}(x'_{q_1}\otimes
x''_{q_2})$, где $x'\in B_{\sigma_1}$, $x''\in B_{\sigma_2}$
\item $(x',x'',x''')=\Phi_{\Omega(\sigma_1,\sigma_2,\sigma_3)}
(x'_{q_1}\otimes x''_{q_2}\otimes x'''_{q_3})$, где $x'\in
B_{\sigma_1}$, $x''\in B_{\sigma_2}$, $x'''\in B_{\sigma_3}$.
\end{itemize}
generate a structure of graph-Frobenius algebra on $B_\star =
\bigoplus_{\sigma\in\Sigma} B_{\sigma}$.
\end{theorem}
\textsl{Proof}. Only the last axiom is not obvious. Let us consider
a film surface
$\Omega\in\Omega(\sigma_1,\sigma_2,\sigma_3,\sigma_4)$, and a
graph-cut between the pairs of vertexes $\sigma_1,\sigma_2$ and
$\sigma_3,\sigma_4$. Then the cut-invariant axiom  and lemma
\ref{l1} give
$\sum_{i,j}(x_1,x_2,b_i^{(1,2|3,4)})F^{ij}_{(1,2|3,4)}(b_j^{(3,4|1,2)},x_3,x_4)
= \Phi_\Omega (x_1,x_2,x_3,x_4)$. Similarly,
$\sum_{i,j}(x_4,x_1,b_i^{(4,1|2,3)})F^{ij}_{(4,1|2,3)}(b_j^{(2,3|4,1)},x_2,x_3)
= \Phi_\Omega (x_1,x_2,x_3,x_4)$

$\Box$

\begin{theorem} \label{t2} Let $B_\star = \bigoplus_{\sigma\in\Sigma}
B_{\sigma}$ be a graph-Frobenius algebra with poli-linear forms
$(.,.)$ and $(.,.,.)$. Then it generates a Film Topological Field
Theory on $\{B_{\sigma}|\sigma\in\Sigma\}$ by means of following
construction. Fix a basis $\{b_i^\sigma\}$ of any vector space
$B_\sigma$, $\sigma\in\Sigma$. Consider the matrix $F^{ij}_\sigma$
that is the inverse matrix for
$F_{ij}^\sigma=(b_i^\sigma,b_j^{\sigma^*})$. Define the linear
functionals on connected film surfaces by
\begin{itemize}
\item $\Phi_{\Omega(\sigma_1,\sigma_2...,\sigma_n)}(x^1_{q_1}\otimes
x^2_{q_2}\otimes ...\otimes
x^n_{q_n})=\sum_{\varsigma_1,\varsigma_2,...,\varsigma_{n-3}\in\Sigma}
(x^1,x^2,b_{i_1}^{\varsigma_1})F_{i_1j_1}^{\varsigma_1}$
$(b_{j_1}^{\varsigma_1^*},x^3_{q_3},b_{i_2}^{\varsigma_2})\\F_{i_2j_2}^{\varsigma_2}
(b_{j_2}^{\varsigma_2^*},x^4_{q_4},b_{i_3}^{\varsigma_3})...
...F_{i_{n-4}j_{n-4}}^{\varsigma_{n-4}}
(b_{j_{n-4}}^{\varsigma_{n-4}^*},x^{n-2}_{q_{n-2}},b_{i_{n-3}}^{\varsigma_{n-3}})
F_{i_{n-3}j_{n-3}}^{\varsigma_{n-3}}
(b_{j_{n-3}}^{\varsigma_{n-3}^*},x^{n-1}_{q_{n-1}},x^n_{q_n})$,\\
where $x^i_{q_i}\in B_{\sigma_i}.$
\end{itemize}
Define the linear functionals on non-connected film surfaces by
multiplicativity axiom.
\end{theorem}
\textsl{Proof}. The  topological invariance follows from the
invariance under cyclic renumbering the vertices of $\Omega$. The
invariance under the renumbering $q_i\mapsto q_j$  $j\equiv
i+1(\mathrm{mod} 2)$ follows from the last axiom for the tree-linear
form. The cut invariance follows directly from the definition of
$\Phi$ if we renumber the vertices marking the cut divide the
vertices $q_1,q_2,...q_k$ and $q_{k+1},q_{k+2},...q_n$.

$\Box$

These two theorems determine the one-to-one correspondence between
Film Topological Field Theories and isomorphic classes of
graph-Frobenius algebras.

\section{Cyclic Foam Topological Field Theories} \label{r3}

\subsection{Cyclic foams} \label{r3.1}

\textit{Cyclic foams} $\Omega$ is defined by a 4-tuple
$(\check{\Omega},\overrightarrow{\partial\Omega},\Delta,\varphi)$,
where
\begin{itemize}
\item $\check{\Omega}=\check{\Omega}(\Omega)$ is a
a compact 2-manifold with a boundary $\partial\check{\Omega}$ that
consists of pairwise non-intersecting circles; moreover, some of
these circles are oriented;
\item $\overrightarrow{\partial\Omega}\subset\partial\check{\Omega}$
is the subset of all oriented circles. The rest circles are
called \textit{free circles};
\item $\Delta=\Delta(\Omega)$ is a regular graph
\item $\varphi:\overrightarrow{\partial\Omega}\to\Delta$ is a \textit{gluing
map}, that is, a homeomorphism on any circle and
$\varphi(\overrightarrow{\partial\Omega})=\Delta$.
\end{itemize}
Here:
\begin{itemize}
\item the cyclic foam $\Omega$ has a cyclic order; this means that the
vertices of any connected component of $\Delta$ have cyclic order
that is agreed with the orientation of
$\varphi(\overrightarrow{\partial\Omega})$;
\item the cyclic foam $\Omega$ is colored; this means that a color
$s(\omega)\in S$ corresponds to each connected component of
$\omega\in\check{\Omega}$ and the colors $s(\omega)$ are pairwise
different for any connected component of $\Omega$;
\item for any connected component $\omega\in\check{\Omega}$
different connected components of
$\partial\omega\cap\overrightarrow{\partial\Omega}$ are mapped
by $\varphi$ to different connected components of $\Delta$;
\item consider a tuple disks $\tilde{\Omega}$ with
$\partial\tilde{\Omega}=\overrightarrow{\partial\Omega}$ and
colors and orientation generated by colors and orientation of
$\overrightarrow{\partial\Omega}$, then the cluing map $\varphi$
generate a film surface $\dot{\Omega}$;
\item a finite set of \textit{marked points} is fixed on
$\Omega$:

(a) marked points from
$\check{\Omega}\setminus\partial\check{\Omega}$ are said to be
\textit{interior} and form a set $\Omega_a$; a point $a\in\Omega_a$
is equipped with a local orientation and the color $s(a)=s(\omega)$
for $a\in\omega\in\check{\Omega}$;

(b) remaining marked points form a set $\Omega_b$ of
\textit{vertices} of $\Omega$; they are all the vertices of
$\Delta$ and the marked points on the free circles; a vertex
graph $\sigma_q$ for the vertex $q\in\Delta$ is defined as
vertex graph for $q\in\dot{\Omega}$; we assume  that each free
circle contains a vertex; the graph of this vertex $q$ is a
segment with an orientation (local orientation of the marked
point) and the color $s(q)=s(\omega)$ for
$q\in\omega\in\check{\Omega}$.
\end{itemize}
Thus the family of cyclic foams contains the family of film surfaces
defined in section \ref{r1} and marked compact 2-manifolds with
boundary considered in \cite{AN}.

We assume 3 types of cuts:
\begin{itemize}
\item a contour-cut, that is a simple closed contour $\gamma\in(\check{\Omega}
\setminus\partial\check{\Omega})$;
\item a segment-cut, a connects without self intersections segment
$\gamma\in\check{\Omega}$ with ends on free contours and without
other intersection with $\partial\check{\Omega}$, that divide
the set of marked points if it divide $\Omega$;
\item a regular graph $\gamma\in\Omega$ that generate a graph-cut on
$\hat\check{\Omega}$.
\end{itemize}

{\ }

Denote by $I_s\in\Sigma$ the isomorphism class of oriented segments
of color $s$. Then $I_s^*=I_s$. Consider families of finite
dimensional vector spaces $\{A_s|s\in S\}$ and
$\{B_{\sigma}|\sigma\in\Sigma\}$. Fix families of tensors
$\{K_s^\otimes\in A_s\otimes A_s|s\in S\}$ and
$\{K_{\sigma}^\otimes\in B_{\sigma}\otimes
B_{\sigma^*}|\sigma\in\Sigma\}$. Fix families of elements
$\{1_{A_s},U_s\in A_s|s\in S\}$ and $\{1_{B_{I_s}}\in B_{I_s}|s\in
S\}$. Fix families of involutions $\{*_s: A_s\rightarrow A_s|s\in
S\}$ and $\{*_{I_s} :B_{I_s}\rightarrow B_{I_s}|s\in S\}$.

Define a functor $\mathcal{V}$ from the category of cyclic foams to
the category of vector spaces. This functor extends the functor on
film surfaces form section \ref{r1}, and the functor on marked
compact 2-manifolds with boundary considered in \cite{AN}.

The  functor $\mathcal{V}$ associates the vector space $V_\Omega=
(\otimes_{p\in\Omega_a}A_{s_p})\otimes(\otimes_{q\in\Omega_b}B_q)$
to any cyclic foam $\Omega$. Here $A_p$ is a copy of $A_{s(p)}$, and
$B_q$ is a copy of $B_{\sigma_q}$. We are going to describe all
morphisms of a monoidal category of cyclic foams and morphisms of
the category of cyclic foams that correspond to it.

\smallskip

(1) {\it Isomorphism.} Let $\phi:\Omega\to \Omega'$ be a
homeomorphism of cyclic foams preserving colors, orientations and
other structures. Define $\mathcal{V}(\phi)=\phi_*: V_{\Omega}\to
V_{\Omega'}$ as a linear operator generated by the bijections
$\phi|_{\Omega_a}:\Omega_b\to\Omega'_a$ and
$\phi|_{\Omega_b}:\Omega_b\to\Omega'_b$.

\smallskip

(2) \textit{Cut.} Let $\Omega$ be a connected film surface and
$\gamma\subset\Omega$ be a cut.

a) Let $\gamma\subset\omega\in\check{\Omega}$ be a non-coorientable
contour-cut. It is presented by a simple coorientable contour
$\gamma'$ on the closure $\overline{\Omega\setminus\gamma}$ of
$\Omega\setminus\gamma$. Contracting  $\gamma'$ to a point $p'$ with
arbitrary local orientation gives the cyclic foam $\Omega'$, where
$V_{\Omega'}=V_{\Omega}\otimes A_{s(\omega)}$. We associate the
morphism $\mathcal{V}(\eta)(x)=\eta_*(x)=x\otimes U_{s(\omega)}$ to
the morphism $\eta:\Omega\to\Omega'$

(b) Let $\gamma\subset\omega\in\check{\Omega}$ be a coorientable
contour-cut. It is presented by simple contours $\gamma_+$ and
$\gamma_-$ on the closure $\overline{\Omega\setminus\gamma}$ of
$\Omega\setminus\gamma$. Contracting  $\gamma_+$ and $\gamma_-$
gives  points $p_+=$ and $p_-$. We assume that their local
orientation are not generated by an orientations of  $\gamma$. Thus
we have a cyclic foam $\Omega'$ and $V_{\Omega'}=V_{\Omega}\otimes
A_{s(\omega)}\otimes A_{s(\omega)}$. We assume the morphism
$\mathcal{V}(\eta)(x)=\eta_*(x)=x\otimes K_{s(\omega)}^\otimes$ to
the morphism $\eta:\Omega\to\Omega'$.

(c) If $\gamma\subset\omega\in\check{\Omega}$ be a segment-cut, then
we define the result of the cutting by $\gamma$ and the value of the
functor on it by analogy with case b) changing  $K_s$ by $K_{I_s}$.

(d) If $\gamma\subset\omega\in\check{\Omega}$ be a graph-cut, then
we define the result of the cutting by $\gamma$ and the value of the
functor on it by analogy with section \ref{r1}.

\smallskip

(3) \textit{Addition of marked point.}

(a)Let us add a not marked point
$p\in\check{\Omega}\setminus\partial\check{\Omega}$ with a local
orientation to the set $\Omega_a$. This operation generates a
morphism $\xi:\Omega\to\Omega'$,  where
$V_{\Omega'}=V_{\Omega}\otimes A_{s(\omega)}$ and
$p\in\omega\in\check{\Omega}$. Associate the morphism
$\mathcal{V}(\xi)(x)=\xi_*(x)=x\otimes 1_{s(\omega)}$ to it.

(b) Similarly, let us add a not marked point
$q\in\overrightarrow{\partial\Omega}
\setminus\partial\check{\Omega}$ with a local orientation to the set
$\Omega_b$. This operation generates a morphism
$\xi:\Omega\to\Omega'$, where $V_{\Omega'}=V_{\Omega}\otimes
B_{I_s(\omega)}$and $p\in\omega\in\check{\Omega}$. Associate the
morphism $\mathcal{V}(\xi)(x)=\xi_*(x)=x\otimes
1_{B_{I{s(\omega)}}}$ to it.

\smallskip

(4) \textit{Change of local orientations of marked points.}

Let $\psi:\Omega\to\Omega'$ be a morphism of change of a local
orientation of a marked point $p\in\Omega_a$ or $q\in\Omega_b$. It
generates an involution $*_{s(p)}: A_p\rightarrow A_p$ or
$*_{I_{s(q)}}: B_q\rightarrow B_q$ and thus the homomorphism
$\mathcal{V}(\phi)=\psi_*: V_{\Omega}\to V_{\Omega'}$.

\smallskip

(5) The tensor product in $\mathcal{S}$ defined by the disjoint
union of surfaces $\Omega'\otimes\Omega''\to \Omega'\coprod\Omega''$
induces the tensor product of vector spaces $\theta_*:
V_{\Omega'}\otimes V_{\Omega''}\to V_{\Omega'\sqcup\Omega''}$.

The functorial properties of $\mathcal{V}$ can be easily verified.

\smallskip

\subsection{Topological Field Theory} \label{r3.2}

Fix families of vector spaces $\{A_s|s\in S\}$ and
$\{B_{\sigma}|\sigma\in\Sigma\}$ and families of tensors, elements
and involutions defining the functor $\mathcal{V}$.

A family of linear forms $\mathcal {F}= \{\Phi_\Omega:V_{\Omega}\to
\mathbb{K}\}$, defined for all cyclic foams $\Omega\in\mathcal{S}$,
is called a \textit{Cyclic Foam  Topological Field Theory} if it
satisfies the following axioms:

\vskip0.6cm $1^\circ$ {\it  Topological invariance.}
$$\Phi_{\Omega'}(\phi_*(x))=\Phi_\Omega(x)$$ for any isomorphism
$\phi:\Omega\to \Omega'$ of cyclic foams.

\vskip 0.6cm $2^\circ$ {\it  Non-degeneracy.}\vskip 0.3cm

Let $\Omega$ (respectively $\Omega^*$) be a sphere with exactly two
marked locally oriented points, where their orientations are induced
by an orientation of the sphere (respectively, there does not exist
an orientation that induces the orientations of the points). Define
the bilinear forms $(.,.)_s, (.,.)_s^*: A_s\times
A_s\rightarrow\mathbb{K}$ by $(x',x'')_s=\Phi_\Omega(x'_{p_1}\otimes
x''_{p_2})$ and $(x',x'')_s^*=\Phi_\Omega^*(x'_{p_1}\otimes
x''_{p_2})$.

Define also the bilinear forms  $(.,.)_{I_s}, (.,.)_{I_s}^*:
B_{I_s}\times B_{I_s}\rightarrow\mathbb{K}$. Their definitions are
similar to the definition of $(.,.)_s, (.,.)_s^*$ after the change
of the sphere by the disk and interior marked points by vertices.
Axiom $2^\circ$ says, that the forms $(.,.)_s, (.,.)_s^*$,
$(.,.)_{I_s}, (.,.)_{I_s}^*$ and the forms $(.,.)_\sigma$ subsection
\ref{r1.2} are non-degenerate.

Let $\Omega$ be a film surface with only two vertices $q_1$, $q_2$.
Then $\sigma_{q_2}=\sigma^*$ if $\sigma_{q_1}=\sigma$. Denote by
$(.,.)_\sigma$ the bilinear form $(.,.)_\sigma: B_\sigma\times
B_{\sigma^*}\rightarrow\mathbb{K}$, where
$(x',x'')_\sigma=\Phi_\Omega(x'_{q_1}\otimes x''_{q_2})$. Axiom
$2^\circ$ asserts that the forms $(.,.)_\sigma$ are non-degenerated
for all $\sigma\in\Sigma$.

\vskip 0.6cm $3^\circ$ {\it Cut invariance.}
$$\Phi_{\Omega'}(\eta_*(x))=\Phi_\Omega(x)$$
for any cut morphism $\eta:\Omega\to \Omega'$ of film surfaces.
\vskip 0.6cm

\vskip 0.6cm $4^\circ$ {\it Invariance under addition of a marked
point.}
$$\Phi_{\Omega'}(\xi_*(x))=\Phi_\Omega(x)$$
for any morphism of addition of a marked point $\xi:\Omega\to
\Omega'$ of film surfaces. \vskip 0.6cm

\vskip 0.6cm $5^\circ$ {\it Invariance under a change of local
orientations.}
$$\Phi_{\Omega'}(\psi_*(x))=\Phi_\Omega(x)$$ for any morphism of
change of the local orientation of a marked point $\psi:\Omega\to
\Omega'$. \vskip 0.6cm

$6^\circ$ {\it Multiplicativity.}

$$\Phi_{\Omega}(\theta^*(x'\otimes x'))= \Phi_{\Omega'}(x')\Phi_{\Omega''}(x'')$$
for $\Omega=\Omega'\cup\Omega''$,  $x'\in V_{\Omega'}$, $x''\in
V_{\Omega''}$ and the morphism of tensor product $\theta:
\Omega'\times\Omega''\to \Omega$.

\vskip 0.6cm

Note that a Topological Field Theory defines the families of
tensors, elements and involutions, defining the functor
$\mathcal{V}$. This follows from lemma \ref{l1} and [\cite{AN2}
lemma 3.1.]

\section{Graph-Cardy-Frobenius algebras}  \label{r4}

\subsection{Definitions} \label{r4.1}

A 3-tuple $(D,l_D,*_D)$ is called an \textit{equipped Frobenius
algebra}(see \cite{AN3}) if $D$ is an associative Frobenius algebra
with unit $1_D$, $l_D: D\rightarrow \mathbb{K}$ is a linear
functional such that the bilinear form$(x_1,x_2)_D=l_D(x_1x_2)$ is
non-degenerate and $*_D: D\rightarrow D$ is an involution such that
$l_D(x^*)=l_D(x)$ and $(x_1x_2)^*=x_2^*x_1^*$ (here and below
$x^*=*(x)$).

Consider a basis $\{d_i|i=1,...,n\}\subset D$, the matrix
$F_{ij}^D=(d_i,d_j)_D$  and the matrix $F^{ij}_D$ inverse to
$F_{ij}^D$. The elements $K_D=F^{ij}_Dd_id_j$ and
$K_D^*=F^{ij}_Dd_id_j^*$ are called the \textit{Casimir} and the
\textit{twisted Casimir} elements, respectively. They don't depend
on the choice of the basis.

We say that a pair of equipped Frobenius algebras
$((A,l_A,*_A),(B,l_B,*_B))$, a homomorphism $\phi: A\rightarrow B$
and an element $U\in A$ form a \textit{Cardy-Frobenius algebra} if
\begin{itemize}
\item $A$ is commutative and the image $\phi(A)$ belongs to the
centre of $B$;
\item $\phi(x^*)=(\phi(x))^*$;
\item $(\phi^*(x),\phi^*(y))_A=\tr W_{x,y}$, where $x,y\in B$,
$(a,\phi^*(b))=(\phi(a),b)_B$, $W\in\End (B)$ and $W(z)=xzy$;
\item $U^2=K_A^*$ and  $\phi(U)=K_B$.
\end{itemize}

It is proved in \cite{AN} that Cardy-Frobenius algebras are in
one-to-one correspondence with Klein Topological Field Theories that
are Topological Field Theories on 2-dimensional manifolds with
boundary. The paper \cite{AN} also contains a complete
classification of semi-simple Cardy-Frobenius algebras.

{\ }

Define now a \textit{graph-Cardy-Frobenius algebra} as a family that
consists of:
\begin{itemize}
\item  a family  of Cardy-Frobenius algebras $\{(A^s,l_A^s,*_A^s)$,
$(B^s,l_B^s,*_B^s),\phi^s,U^s|s\in S\}$
\item a graph-Frobenius algebra $B_\star$ with a bilinear form
$(.,.)_B: B_\star\times B_\star\rightarrow\mathbb{K}$, and a
three-linear form $(.,.,.)_B: B_\star\times B_\star\times
B_\star\rightarrow\mathbb{K}$
\item a family  of homomorphisms $\{\phi_\sigma^s: A^s\rightarrow
\End(B_\sigma)|s\in S,\sigma\in \Sigma\}$, where $\phi_\sigma^s=0$
if $s$ is not the color of an edge of $\sigma$.
\end{itemize}

Here:
\begin{itemize}
\item  $B^s$ coincides with
$B_{I_s}\subset B_\star$ and $\phi_{I^s}^s(a)(b)=\phi^s(a)b$ for
$a\in A^s$, $b\in B^s$;
\item
$(\phi_{\sigma_1}^s(a)(x_1),x_2)_B=(x_1,\phi_{\sigma_2}^s(a)(x_2))_B$,
where $a\in A^s$, $x_i\in B_{\sigma_i}$;
\item $(\phi_{\sigma_1}^s(a)(x_1),x_2,x_3)_B=(x_1,\phi_{\sigma_2}^s(a)(x_2),x_3)_B=
(x_1,x_2,\phi_{\sigma_3}^s(a)(x_3))_B$, where $a\in A^s$, $x_i\in
B_{\sigma_i}$.
\end{itemize}

\subsection{One-to-one correspondence} \label{r4.2}

Let  $\mathcal {F}= \{\Phi_\Omega:V_{\Omega}\to \mathbb{K}\}$ be a
Cyclic Foam Topological Field Theory. Then its restriction to
2-dimensional manifolds with boundary forms a Klein Topological
Field Theory $\mathcal {F}_K$ and, therefore, a family of
Cardy-Frobenius algebras $\{((A^s,l_A^s,*_{A^s}),
(B^s,l_B^s,*_{B^s}),\phi^s,U^s)|s\in S\}_{\mathcal {F}}$. The
restriction to film surfaces forms Film Topological Field Theory
$\mathcal {F}_N$ and, therefore, a  graph-Frobenius algebra
$(B_\star, (.,.)_B, (.,.,.)_B)_{\mathcal {F}}$.

Let us define the homomorphisms $\{\phi_\sigma^s: A^s\rightarrow
\End(B_\sigma)|s\in S,\sigma\in \Sigma\}_{\mathcal {F}}$. Let
$\Omega$ be a cyclic foam with two vertices $q_1$, $q_2$ and one
interior marked point $p$. Put $\sigma=\sigma_{q_1}$. The functional
$\Phi_\Omega$ generates the homomorphism $\phi_\sigma^s$ of $A_s$ to
the space $E$ of linear functionals on $B_\sigma\otimes
B_{\sigma^*}$. The bilinear form $(.,.)_\sigma$ generates the
isomorphism between $B_{\sigma^*}$ and the space of linear
functionals on $B_\sigma$. Thus we can identify $E$ with
$\Hom(B_\sigma,B_\sigma)$.

\begin{theorem}\label{t3} Let $\mathcal {F}=
\{\Phi_\Omega:V_{\Omega}\to \mathbb{K}\}$ be a cyclic Foam
Topological Field Theory on families of vector spaces $\{A_s|s\in
S\}$, $\{B_{\sigma}|\sigma\in\Sigma\}$. Then the Cardy-Frobenius
algebras $\{((A^s,l_A^s,*_{A^s}),
(B^s,l_B^s,*_{B^s}),\phi^s,U^s)|s\in S\}_{\mathcal {F}}$, the
graph-Frobenius algebra $(B_\star, (.,.)_B, (.,.,.)_B)_{\mathcal
{F}}$ and the homomorphisms $\{\phi_\sigma^s: A^s\rightarrow
\End(B_\sigma)|s\in S,\sigma\in \Sigma\}_{\mathcal {F}}$ form a
graph-Cardy-Frobenius algebra.
\end{theorem}
\textsl{Proof}. The properties $B^s=B_{I_s}$, and
$\phi_{I^s}^s(a)(b)=\phi^s(a)b$ follow from the corresponding axiom.
Let us prove that $(\phi_{\sigma_1}^s(a)(x_1),x_2)_B=
(x_1,\phi_{\sigma_2}^s(a)(x_2))_B$. Consider a cyclic foam $\Omega$
that is the film surface $\Omega(\sigma_1,\sigma_2)$ with an
interior marked point $p$ of color $s$. Then the cut axiom gives
$(\phi_{\sigma_1}^s(a)(x_1),x_2)_B=\Phi_\Omega(a\otimes x_1\otimes
x_2)$ and $(x_1,\phi_{\sigma_2}^s(a)(x_2))_B=\Phi_\Omega(a\otimes
x_1\otimes x_2)$. A prove of the identities
$(\phi_{\sigma_1}^s(a)(x_1),x_2,x_3)_B=
(x_1,\phi_{\sigma_2}^s(a)(x_2),x_3)_B=
(x_1,x_2,\phi_{\sigma_3}^s(a)(x_3))_B$ is similar.

$\Box$
\begin{theorem}\label{t4} The correspondence from theorem \ref{t3}
generates a one-to-one correspondence between Cyclic Foam
Topological Field Theories and isomorphism classes of
graph-Cardy-Frobenius algebras.
\end{theorem}
\textsl{Proof}. Let $(\{((A^s,l_A^s,*_{A^s}),
(B^s,l_B^s,*_{B^s}),\phi^s,U^s)|s\in S\}$, $(B_\star, (.,.)_B,
(.,.,.)_B)$, $\{\phi_\sigma^s: A^s\rightarrow \End(B_\sigma)|s\in
S,\sigma\in \Sigma\})$ be a graph-Cardy-Frobenius algebra. Let us
construct a Cyclic Foam Topological Field Theory $\mathcal {F}=
\{\Phi_\Omega:V_{\Omega}\to \mathbb{K}\}$ that generates it.
According to the cut axiom and the axiom
$\phi_{I^s}^s(a)(b)=\phi^s(a)b$, the Theory $\mathcal {F}$ is
defined by its restrictions to:
\begin{itemize}
\item  2-dimensional manifolds with boundary and arbitrary number of
marked points;
\item film surfaces without interior marked points;
\item film surfaces with two vertices and one interior marked point.
\end{itemize}
According to \cite{AN}, Topological Field Theories on 2-dimensional
manifolds with boundary and arbitrary number of marked points are in
one-to-one correspondence with isomorphism classes of
Cardy-Frobenius algebras $(\{((A^s,l_A^s,*_{A^s}),
(B^s,l_B^s,*_{B^s}),\phi^s,U^s)|s\in S\}$. According to \ref{t2},
Topological Field Theories on film surfaces without interior marked
points  are in one-to-one correspondence with isomorphism classes of
graph-Frobenius algebras $(B_\star,(.,.)_B, (.,.,.)_B)$. Define the
value of $\mathcal {F}$ on surfaces $\Omega$ with two vertices and
one interior marked point by $\Phi_\Omega(a\otimes x_1\otimes
x_2)=(\phi_\sigma^s(a)(x_1), x_2)_B$. The properties
$\phi_{I^s}^s(a)(b)=\phi^s(a)b$,
$(\phi_{\sigma_1}^s(x_1),x_2)_B=(x_1,\phi_{\sigma_2}^s(x_2))_B$ and
$(\phi_{\sigma_1}^s(a)(x_1),x_2,x_3)_B=(x_1,\phi_{\sigma_2}^s(a)(x_2),x_3)_B=
(x_1,x_2,\phi_{\sigma_3}^s(a)(x_3))_B$ guarantee that the values of
$\mathcal {F}$ satisfy the axiom of Cyclic Foam Topological Field
Theory.

$\Box$

\begin{note} The category of cyclic foams contain the
subcategory of \textsf{oriented foams} $\Omega=
(\check{\Omega},\overrightarrow{\partial\Omega},\Delta,\varphi)$,
where the orientation of $\overrightarrow{\partial\Omega}$ generate
orientations of edges of $\Delta$. Our constructions make possible
to define the Topological Field Theories for oriented foams and to
prove that these Topological Field Theories one-to-one correspond to
the analog of graph-Cardy-Frobenius algebras where arbitrary colored
graphes change to bipartite colored graphs.
\end{note}

\section{Examples of Cyclic Foam  Topological Field Theories} \label{r5}

In this section we construct an example of Cyclic Foam  Topological
Field Theory. Its restriction to 2-dimensional manifolds with
boundary is the Klein  Topological Field Theory of regular covering,
constructed in \cite{AN3}.

Associate a group $G_s$ of a set $X_s$ and an action of $G_s$ on
$X_s$ to any color. Consider the vector space $A_s$ which is the
center of the group algebra of $G_s$. Associate
$X_{\tilde{S}}=\times_{s\in S}X_s$ to any finite subset
$\tilde{S}\subset S$. The actions of $G_s$ on $X_s$ generate the
action of $G=\bigoplus_{s\in S}G_s$ on $X_{\tilde{S}}$.

Let $L=L(\tilde{\sigma})$ be the set of edges of a colored graph
$\tilde{\sigma}$. Let $s(l)$ be the color of $l\in L$. Denote by
$\tilde{\sigma}^X$ the set of all maps $\psi: L\rightarrow
X_{s(L)}\times X_{s(L)}$, where $\psi(l)\in X_{s(l)}\times
X_{s(l)}$. Define the action of $G$ on $\tilde{\sigma}^X$ by
$g(\psi(l))=g(x')\times g(x'')$ for $\psi(l)=x'\times x''$. Let
$\tilde{\sigma}^{X_G}$ be the set of orbits of this action.

A pair $(\tilde{\sigma},\psi_G)$, where $\tilde{\sigma}$ is a
colored graph and $\psi_G\in\tilde{\sigma}^{X_G}$, is called an
\textit{equipped colored graph}, or a colored graph with
\textit{equipment} $\psi_G$. An isomorphism
$\varphi:\tilde{\sigma}_1\rightarrow\tilde{\sigma}_2$ of colored
graphs is called an \textit{isomorphism }of the equipped colored
graphs $(\tilde{\sigma}_i,\psi_G^i)$ if it takes $\psi_G^2$ to
$\psi_G^1$. Denote by $|\Aut(\sigma,\psi_G)|$ the order of the group
$\{g\in G| g\psi =\psi\}$, where
$\psi\in\psi_G\in\tilde{\sigma}^{X_G}$ and
$\tilde{\sigma}\in\sigma\in\Sigma$.

Consider the set $E_{\tilde{\sigma}}$  of equipments of a colored
graph $\tilde{\sigma}$. Isomorphisms of colored graphs generate the
canonical bijections between the corresponding sets
$E_{\tilde{\sigma}}$. Thus we can associate the set $E_\sigma$ to
any $\sigma\in\Sigma$. Let $B_\sigma$ be a vector space generated by
$E_\sigma$. Denote by $*: B_\sigma\rightarrow B_{\sigma^*}$ the
involution generated by changing orientation of $\sigma$ and
changing of the components of $X_s\times X_s$ .

{\ }

Construct a  Cyclic Foam  Topological Field Theory with families of
vector spaces $\{A_s|s\in S\}$, $\{B_\sigma|\sigma\in\Sigma\}$, by
its restriction to
\begin{itemize}
\item  2-dimensional manifolds with boundary and arbitrary number of
marked points;
\item film surfaces without interior marked points;
\item film surfaces with two vertices and one interior marked point.
\end{itemize}

We start with its description of $\mathcal{F}$ on film surfaces.
Consider an additional structure on film surfaces. Let $V=V(\Omega)$
be the set of edges of a film surface $\Omega$. Denote by $S(v)$ the
set of the colors of the disks that are incident to $v\in V$.
Consider the set $\Omega^X$ of maps $\psi: V\rightarrow
\bigcup_{v\in V}X_{S(v)}$, where $\psi(v)\in X_{S(v)}$. Denote by
$\Omega^X$ the set of orbits for the action of $G$ on $\Omega^X$.

A pair $(\Omega,\psi_G)$, where $\Omega$ is a colored graph and
$\psi_G\in\Omega^{X_G}$, is called an \textit{equipped film surface}
or film surface with \textit{equipment} $\psi_G$. An isomorphism
$\varphi:\Omega_1\rightarrow\Omega_2$ of film surfaces is called an
\textit{equivalence} of the equipment film surfaces
$(\Omega_i,\psi_G^i)$ if it takes $\psi_G^2$ to $\psi_G^1$. Denote
by $|\Aut(\Omega,\psi_G)|$ the order of the group $\{g\in G| g\psi
=\psi\}$, where $\psi\in\psi_G\in\Omega^{X_G}$. An equipment $\psi$
of the film surface $\Omega$ generates an equipment $\psi_\sigma$
for the graph $\sigma$ of any vertex of $\Omega$. We assume that
$\psi_\sigma(l)=(x^1_{s(l)},x^2_{s(l)})$, where $l\in L(\sigma)$ is
an oriented edge from $v^1\in V(\Omega)$ to $v^2\in V(\Omega)$ and
$\psi(v^i)=\times_{s\in S(v^i)}x^i_s$.

We say that a connected equipped film surface $(\Omega,\psi_G)$ is a
\textit{compatible surface} for equipped colored graphes
$\varsigma_1,\varsigma_2,...,\varsigma_n$ if these graphes are the
vertex graphs of $\Omega$ and the numeration the graphs $\sigma_i$
generates the cyclic order of vertexes of film surface $\Omega$.
Denote by $\Psi(\varsigma_1,\varsigma_2,...,\varsigma_n)$ the set of
all isomorphism classes of compatible surfaces for
$\varsigma_1,\varsigma_2,...,\varsigma_n$.

Define a set of linear functionals $\mathcal {F}_N=
\{\Phi_\Omega:V_{\Omega}\to \mathbb{K}\}$ on connected equipped film
surfaces by $ \Phi_\Omega(\varsigma_1\otimes
\varsigma_2\otimes...\otimes\varsigma_n)=\sum_{\Psi\in
\Psi(\varsigma_1,
\varsigma_2,...,\varsigma_n)}\frac{1}{|\Aut(\Psi)|}$.
\begin{lemma}\label{l2} The set $\mathcal {F}_N=
\{\Phi_\Omega:V_{\Omega}\to \mathbb{K}\}$ generates a Film
Topological Field Theory.
\end{lemma}

\textsl{Proof} It follows from our definition that
$(\varsigma_1,\varsigma_2^*)_\sigma=
\frac{\delta_{\varsigma_1,\varsigma_2^*}} {|\Aut(\varsigma_1)|}$ for
$\varsigma_1,\varsigma_2\in E_\sigma$, and thus
$K_\sigma=\sum_{\varsigma\in
E_\sigma}|\Aut(\varsigma)|\varsigma\otimes\varsigma^*$. Let us prove
the cut invariance.

Consider an equipment film surface
$\Psi\in\Psi(\varsigma_1,\varsigma_2,...,\varsigma_n)=
(\Omega(\sigma_1,\sigma_2,...,\sigma_n),\psi_G)$. Let $\eta$ be the
cut morphism by graph-cut
$\gamma\subset\Omega(\sigma_1,\sigma_2,...,\sigma_n)$. It associates
the pair of film surfaces
$\Omega(\sigma_1,\sigma_2,...,\sigma_k,\sigma')$,
$\Omega(\sigma'',\sigma_{k+1},\sigma_2,...,\sigma_n)$ to the film
surface $\Omega(\sigma_1,\sigma_2,...,\sigma_n)$. Any equipment of
$\Omega(\sigma_1,\sigma_2,...,\sigma_n)$ generates equipments of
$\Omega(\sigma_1,\sigma_2,...,\sigma_k,\sigma')$,
$\Omega(\sigma'',\sigma_{k+1},\sigma_2,...,\sigma_n)$. Thus we
receive equipment film surfaces $\Psi'$ и $\Psi''$ and an equipment
of $\sigma'$.

Fix an equipped colored graph $\varsigma'=(\sigma',\psi'_G)$  and
consider the set of equipped film surfaces
$\Psi_{\varsigma'}(\varsigma_1,\varsigma_2,...,\varsigma_n)\subset
\Psi(\varsigma_1,\varsigma_2,...,\varsigma_n)$ that generate the
equipment $\psi'_G$ on $\sigma'$. Then
$\sum_{\Psi\in\Psi_{\varsigma'}(\varsigma_1,\varsigma_2,...,\varsigma_n)}
\frac{1}{|\Aut(\Psi)|}=$
$\frac{|\Aut(\varsigma)|}{|\Aut(\Psi')||\Aut(\Psi'')|}$. Summations
over all equipments $\psi'_G$ of $\sigma'$ gives
$\Phi_\Omega(\varsigma_1\otimes
\varsigma_2\otimes...\otimes\varsigma_n)=$
$\sum_{\Psi\in\Psi_\sigma(\varsigma_1,\varsigma_2,...,\varsigma_n)}
\frac{1}{|\Aut(\Psi)|}=$ $\sum_{\varsigma\in E_\sigma}
\frac{|\Aut(\varsigma)|}{|\Aut(\Psi_+)||\Aut(\Psi_-)|}=$
$\Phi_{\Omega'}\eta_*(\varsigma_1\otimes
\varsigma_2\otimes...\otimes\varsigma_n)$.

$\Box$

It follows from the previous section that $\mathcal {F}_N$ is
generated and is defined by the graph-Frobenius algebra $B_*$. This
algebra has the basis $E=\bigcup_{\sigma\in\Sigma}E_\sigma$ and is
defined by polilinear forms
\begin{itemize}
\item $(\varsigma_1,\varsigma_2)_B=
\sum_{\Psi\in\Psi(\varsigma_1,\varsigma_2)}
\frac{1}{|\Aut(\Psi)|}=\frac{\delta_{\varsigma_1,\varsigma_2^*}}{|\Aut(\varsigma_1)|}$
for $\varsigma_1,\varsigma_2\in E$;
\item $(\varsigma_1,\varsigma_2,\varsigma_3)_B=
\sum_{\Psi\in\Psi(\varsigma_1,
\varsigma_2,\varsigma_3)}\frac{1}{|\Aut(\Psi)|}$ for
$\varsigma_1,\varsigma_2,\varsigma_3\in E$.
\end{itemize}

{\ }

Define now the action of the group algebra $A_s$ on $B_\sigma$. It
is identical if there are no edges of color $s$ between the edges of
$\sigma$. Let $s(l)=s$, $\psi\in\psi_G\in\sigma^{X_G}$,
$\psi(l)=(x',x'')$ and $a=\Sigma_{g\in G}\lambda_g g$. Then we
assume that $\phi_\sigma^s(a)(\psi)=(\Sigma_{g\in G}\lambda_g
gx',x'')$ on $l$ and $\phi_\sigma^s(a)(\psi)=\psi$ on the other
edges of $\sigma$. The function $\phi_\sigma^s(a)(\psi)$ depends
only on its orbit $\psi_G$, and thus generate the linear operator
$\phi_\sigma^s(a): B_\sigma\rightarrow B_\sigma$.

Define now the system of linear operators $\mathcal {F}_C=
\{\Phi_\Omega:V_{\Omega}\to \mathbb{K}\}$ on cyclic foams with two
vertices by $\Phi_\Omega(a^1\otimes...\otimes a^r\otimes x_1\otimes
x_2)= (\phi_\sigma^{s(a^1)}(a^1)...\phi_\sigma^{s(a^r)}(a^r)(x_1),
x_2)_B$.

Define a system of linear operators $\mathcal {F}_s=
\{\Phi_\Omega:V_{\Omega}\to \mathbb{K}\}$ on 2-dimensional manifolds
with boundary of color $s$. We set it to be the Klein Topological
Field Theory of $G_s$-regular covering with trivial stationary
subgroup from  \cite{AN3}.

\begin{theorem}\label{t5} There exists a unique Cyclic Foam  Topological Field
Theory $\mathcal {F}$ with the restrictions $\mathcal {F}_N$,
$\mathcal {F}_C$ и $\mathcal {F}_s$.
\end{theorem}

\textsl{Proof} It follows from lemma \ref{l2} and \cite{AN3} that
the families $\mathcal {F}_N$ and $\mathcal {F}_s$ satisfy the
axioms of Cyclic Foam  Topological Field Theory. By our definitions,
the value of $\mathcal {F}_C$ on $\Omega$ is equal to the product of
the values $\mathcal {F}_s$ on the disks that form $\Omega$. Thus
$\mathcal {F}_C$ also satisfies the axioms of Cyclic Foam
Topological Field Theory. Moreover, the families $\mathcal {F}_N$,
$\mathcal {F}_C$, and $\mathcal {F}_s$ coincide on common areas of
definition.

To define $\mathcal {F}$ on an arbitrary cyclic foam one can use the
cut axiom and cut the surface into 2-dimensional manifolds with
boundary, film surfaces without interior marked points and surfaces
with two vertices. The result does not depend on the cut system
because any two such systems are different only on 2-dimensional
manifolds with boundary, film surfaces without interior marked
points or surfaces with two vertices.

$\Box$


\begin{thebibliography}{References}

\bibitem{AN} Alexeevski A., Natanzon S., Noncommutative two-dimensional
topological field theories and Hurwitz numbers for real algebraic
curves. Selecta Math., New ser. v.12,n.3, 2006, p. 307-377 (arXiv:
math.GT/0202164).

\bibitem{AN1} Alexeevski A., Natanzon S., Algebra of Hurwitz numbers for
seamed surfaces, Russian Math.Surveys, 61 (4) (2006), 767-769

\bibitem{AN2} Alexeevski A., Natanzon S., Algebra of bipartite graphs and
Hurwitz numbers of seamed surfaces. Accepted to Math.Russian
Izvestiya

\bibitem{AN3} Alexeevski A., Natanzon S., Hurwitz numbers for
regular coverings of surfaces by seamed surfaces and Cardy-Frobenius
algebras of finite groups, arXiv: math/07093601

\bibitem{At} Atiyah M., Topological Quantum Field Theories, Inst. Hautes
Etudes Sci. Publ. Math., 68 (1988), 175-186.

\bibitem{B} Baez J.C., An introdaction to Spin Foam Models of BF
Theory and Quantum Gravity, arXiv:gr-qc/9905087

\bibitem{CNP} Costa A.F., Natanzon S.M., Posto A.M., Counting the regular
coverings of surfaces using the center of a group algebra, European
Journal of combinatorics, 27 (2006), 228-234.

\bibitem{D2} Dijkgraaf  R., Geometrical Approach to Two-Dimensional
Conformal Field Theory, Ph.D.Thesis (Utrecht, 1989)

\bibitem{D1} Dijkgraaf  R., Mirror symmetry and elliptic curves, The moduli
spaces of curves, Progress in Math., 129 (1995), 149-163,
Brikh\"auser.

\bibitem{Du} B.Dubrovin, Geometry of 2D topological field theories In: LNM,
1620 (1996), 120-348.

\bibitem{F} Faith C., Algebra II Ring theory, Springer-Verlag, 1976

\bibitem{H} Hurwitz A., \"Uber Riemann'sche Fl\"achen mit gegeben
Verzweigungspunkten, Math., Ann., Bn.39 (1891), 1-61.

\bibitem{K} Khovanov M. l(3) link homology I,
arXiv: math.QA/0304375.

\bibitem{KR1} Khovanov M., Rozansky L., Matrix factorizations and link
homology hep-th/0401268.

\bibitem{KR2} Khovanov M., Rozansky L., Topological Landau-Ginzburg models
on a world-sheet foam., arXiv: hep-th/0404189.

\bibitem{Lauda} Lauda A.D., Pfeiffer H. Open-closed strings: Two-dimensional extended TQFTs
and Frobenius algebras. arXiv:math.AT/0510664.

\bibitem{Laz} Lazaroiu C.I., On the structure of open-closed topological
field theory in two-dimensions, Nucl. Phys. B 603 (2001), 497-530.

\bibitem{Moore} Moore G., Some comments on branes, G-flux, and K-theory,
Int.J.Mod.Phys.A 16,936(2001), arXiv:hep-th/0012007

\bibitem{MS} Moore G., Segal G, D-branes and K-theory in 2D topological field theory,
arXiv:hep-th/0609042 (2006)

\bibitem{O} Oriti D., Spin foam models of quantum spacetime
arXiv:gr-qc/0311066 (2003)

\bibitem{PT} Porter.T, Turaev.V., Formal homotopy quantum  field
theories,I: formal maps and crossed C-algebras arXiv:hep-th/0512032.

\bibitem{R} Rozansky L., Topological A-models on seamed Riemann surfaces.
arXiv: hep-th/0305205.

\bibitem{Se} Segal G.B., Two dimensional conformal field theory and
modular functor. In: Swansea Proceedings,Mathematical Physics, 1988,
22-37.

\bibitem{TT} Turaev V., Turer P., Unoriented topological quantum  field
theory and link homology, Algeb. Geom. Topol.,6(2006).

\bibitem{W} Witten.E, Quantum Field theory.
Commun.Math.Phys.117(1988),353-386


\end{thebibliography}
\end{document}